\RequirePackage{fix-cm}
\documentclass[twocolumn]{svjour3}          
\smartqed  
\usepackage{graphicx}
\usepackage{enumerate}
\usepackage[utf8]{inputenc}
\usepackage{amsmath,amssymb}
\usepackage{enumitem}
\usepackage{times}
\usepackage{graphicx} 
\usepackage{eucal}
\usepackage{amsfonts}
\usepackage{xcolor}
\usepackage{comment}
\usepackage{float}
\usepackage{enumitem}
\usepackage{hyperref}
\usepackage{xspace}
\usepackage{caption}
\usepackage{subcaption}
\usepackage{amssymb}
\usepackage{latexsym}
\usepackage{appendix}
\usepackage{pdfpages}
\journalname{CGI2024} 

\newcommand{\newtext}[1]{{#1}\xspace}

\newcommand{\ie}{i.e.,\xspace}
\newcommand{\etal}{et al.\xspace}
\newcommand{\myparagraph}[1]{\noindent\textbf{\\
#1}}
\usepackage [english]{babel}
\usepackage [autostyle, english = american]{csquotes}
\MakeOuterQuote{"}

\begin{document}

\title{Jacobi Set Simplification for Tracking Topological Features in Time-Varying Scalar Fields}
\subtitle{}
\author{Dhruv Meduri \and Mohit Sharma \and Vijay Natarajan}
\institute{Dhruv Meduri (Corresponding author) \at University of Utah, USA.  \email{u1471195@utah.edu} \and Mohit Sharma and Vijay Natarajan\at Indian Institute of Science, Bangalore, India\\ \email{\{mohitsharma,vijayn\}@iisc.ac.in}}
\date{ }

\maketitle

\begin{abstract}
The Jacobi set of a bivariate scalar field is the set of points where the gradients of the two constituent scalar fields align with each other. It captures the regions of topological changes in the bivariate field. The Jacobi set is a bivariate analog of critical points, and may correspond to features of interest. In the specific case of time-varying fields and when one of the scalar fields is time, the Jacobi set corresponds to temporal tracks of critical points, and serves as a feature-tracking graph. The Jacobi set of a bivariate field or a time-varying scalar field is complex, resulting in cluttered visualizations that are difficult to analyze. This paper addresses the problem of Jacobi set simplification. Specifically, we use the time-varying scalar field scenario to introduce a method that computes a reduced Jacobi set. The method is based on a stability measure called robustness that was originally developed for vector fields and helps capture the structural stability of critical points. We also present a mathematical analysis for the method, and describe an implementation for 2D time-varying scalar fields. Applications to both synthetic and real-world datasets demonstrate the effectiveness of the method for tracking features.

\end{abstract}

\section{Introduction}
\label{sec:introduction}
The relationship between scalar fields defined over a spatial domain of interest corresponds to interesting phenomena or captures application-specific features. For example, a water mass in oceanography is defined by the temperature-salinity (T-S) curve~\cite{helland1916nogenTSDiag},  fields such as pressure and wind speed are used to identify structural characteristics of a hurricane~\cite{doleisch2004IsabelScivis}. This motivates the study of techniques for visualizing multifield or multivariate data~\cite{hansen2014scientific}.  Various structures have been introduced to analyze the relationship between multiple fields and specifically for bivariate fields, such as continuous scatterplot~\cite{Bachthaler2008CSP}, fiber surface~\cite{carr2015fiber}, Reeb space~\cite{Edelsbrunner2008ReebSpace}, Pareto set~\cite{huettenberger2014decomposition}, and Jacobi set~\cite{edelsbrunner2004jacobi}. We restrict our attention to the Jacobi set, a topological descriptor that is a generalization of the notion of critical points to the multivariate setting. For a bivariate field, the Jacobi set is defined as the set of points in the domain where the gradients of the two fields align with each other.

This paper focuses on time-varying scalar fields that may be represented as a bivariate field over space-time, where the time function is one of the two scalar fields. The Jacobi set of a time-varying scalar field represents the tracks of critical points of the scalar field over time. Critical points may correspond to interesting features, such as vortices in fluid flows or eddies in the ocean, and hence their tracks are of interest. The Jacobi set is, in this case, a feature-tracking graph. The Jacobi set computed from piecewise linear (PL) scalar functions is sensitive to noise, dense, and cluttered, making it difficult to analyze its structure. We propose a method that computes a simplified Jacobi set, resulting in a smaller number of clutter-free tracks while retaining important features that are represented by the original Jacobi set. The intuition from time-varying fields is central to the development of our simplification method and we consider this as a first step towards the development of methods that are applicable to generic bivariate fields.

We present a Jacobi set simplification method that considers both temporal and spatial coherence of features. It formulates critical points of a scalar field as zeros of the corresponding gradient field thereby enabling a bottom-up simplification directed by a stability measure called robustness that was introduced for vector fields~\cite{Primoz2015RobustnessSimplification4VF}. Clusters of critical points that are determined by the robustness parameter are tracked over time. The track of a cluster of critical points replaces the tracks of all critical points within the cluster thereby reducing the number of tracks while retaining their spatial and temporal structure. The resulting collection of tracks constitutes the simplified Jacobi set. Key contributions of this paper include
\begin{enumerate}
    \item A robustness-based Jacobi set simplification algorithm for time-varying scalar fields,
    \item A mathematical analysis supporting the algorithm that guarantees the existence of a corresponding simplified vector field, and
    \item An implementation of the method for 2D time-varying scalar fields.
\end{enumerate}
The effectiveness of the method is demonstrated via experiments on a synthetic dataset and multiple fluid flow datasets. 

\vspace{-0.1in}

\section{Related work}
This paper focuses on simplifying the Jacobi set for time-varying scalar fields, resulting in reduced clutter in the visualizations of feature tracks. This section summarizes previous work on structures, including the Jacobi set, that have been employed for bivariate field visualization, different applications of the Jacobi set, methods for simplification, and feature tracking.

\myparagraph{Multivariate field visualization.}
Different concepts and structures have been introduced to study the relationship between individual fields in a multivariate dataset. Continuous scatterplot~(CSP)~\cite{Bachthaler2008CSP}, a generalization of scatterplot to continuous fields, provides a dense visualization in range space. This visual representation helps identify regions of interest that may in turn be mapped to the spatial domain using fiber surfaces~\cite{carr2015fiber}. CSPs and fiber surfaces have been used for studying bivariate and multivariate fields from multiple application domains~\cite{Blecha2019Nuclear,Raith2019Tensor,sharma2021segmentation,Sharma2023CSPOperators}. Nagaraj and Natarajan~\cite{Nagaraj2011isoExtractMultifield} introduce a variation density function that captures the relationship between multiple scalar fields over isosurfaces of a given scalar field. This profile helps identify interesting isovalues of scalar fields in multivariate scenarios. Chattopadhyay~\etal~\cite{Chattopadhyay2014Jacobist} introduce the Jacobi structures as the critical features
in the Reeb space~\cite{Edelsbrunner2008ReebSpace}, establishing a relationship between the Jacobi set and singular fibers. The Reeb space is the bivariate analog of the Reeb graph. A fiber component, the bivariate counterpart of an isocontour component, maps to a point in the Reeb space. The joint contour net~\cite{Carr2014JCN} is a discrete representation of the Reeb space. The representation and computation of Reeb space is challenging, and hence it is often not computed explicitly in practice. The Jacobi set is a subset, specifically the image of singular fibers, of the Reeb space and is amenable to fast computation. The notion of Pareto set~\cite{huettenberger2014decomposition} helps identify consensus regions using Pareto optimality and dominance in multifields.

\myparagraph{Jacobi set applications.} 
The Jacobi set has been applied for many tasks in image analysis and visualization. Edelsbrunner~\etal~\cite{Edelsbrunner2004LnGcompJS} introduce local and global comparison measures based on the Jacobi set to support comparative visualization. In image processing, Norgard and Bremer~\cite{Norgard2013ridgeValley} utilize the Jacobi set to extract image ridges, proposing a combinatorial algorithm to create a ridge-valley graph that captures information about all ridges in a given image. 
The Jacobi set has been used in geophysics for estimating relationships between geopotential height and total ozone column~\cite{Artamonova2017physicsJS} and in forestry for automatic tree ring detection~\cite{Makela2020AutomaticTR}. Tierny and Carr use the Jacobi fiber surfaces, the preimage of a Jacobi edge, to compute the Reeb space of a bivariate field~\cite{tierny2016jacobi}. Sharma and Natarajan~\cite{sharma2022FF} propose an output-sensitive algorithm for fiber surface computation that further facilitates the computation of individual fiber surface components in the vicinity of a Jacobi edge selected by the user. The performance of the algorithm depends on the size of the Jacobi set, and it will benefit from a simplified Jacobi set.

\myparagraph{Jacobi set simplification.}
 In practical scenarios, the Jacobi set consists of a complex network of edges, including noisy artifacts that result from the PL representation. Analyzing it can be challenging and will benefit from a simplification method. Simplifying the bivariate field will, indirectly, result in a simpler Jacobi set. The bivariate field may be simplified by employing topological simplification directed by persistence on the \newtext{individual scalar fields~\cite{ELZ02,bremer2007topological,lukasczyk2020localized}.} However, in a time-varying scenario, such an approach may introduce temporal inconsistencies. Post simplification, a critical point may no longer have a correspondence in adjacent time steps, breaking the temporal continuity.
 The Jacobi set is not a stable structure. Small changes in the underlying bivariate field can produce significant changes in the Jacobi set. This instability poses challenges for such indirect simplification, especially in time-varying fields.
 
Direct simplification methods aim to directly remove noisy edges or sections of Jacobi curves that are deemed to be less significant. These approaches are directed by a metric that measures the significance of edges of the Jacobi set. Suthambhara and Natarajan~\cite{Suthambhara2011jacobisetsimplification} propose a method to reduce the number of connected components in the Jacobi set by posing the simplification problem as the extraction of level sets and offset contours. Bhatia~\etal~\cite{bhatia2015localjacobisetsimplification} generalize the concept of critical point cancellation to Jacobi sets and simplify 'Jacobi regions'. The method identifies less significant portions of the Jacobi set and removes them by modifying the individual scalar fields locally. Recently developed methods~\cite{klotzl2022local,klotzl2022reduced} compute a clutter-free visual representation of the Jacobi set by considering a bilinear interpolant. The method removes zig-zag artifacts while maintaining the topological structure.

Our proposed method may be classified under the direct simplification category. It reduces the number of Jacobi edges by clustering edges that potentially represent a common feature for a given approximation threshold. It directly computes a simplified version of the Jacobi set as opposed to simplifying the noisy Jacobi set.

\myparagraph{Feature tracking.} 
In a time-varying scalar field, the Jacobi edges represent tracks of critical points. Since critical points are good representatives of topological features, critical point tracking is a preferred approach for topological feature tracking. Feature tracking has been studied widely, but it has not been formulated using the Jacobi set except for a few instances~\cite{edelsbrunner2008time}. Reininghaus~\etal~\cite{reininghaus2011efficient} use ideas from discrete Morse theory and combinatorial feature flow fields to track critical points. They introduce the notion of integrated persistence, combining spatial persistence of critical points with their temporal evolution to determine the significance of a particular track as opposed to direct elimination of short-lived features. Other works~\cite{lukasczyk2017nested,saikia2017global,Lukasczyk2019DNTG,Yan2021Survey} focus on tracking regions of interest, such as predefined connected components instead of critical points. \newtext{Feature tracking has also been studied for vector fields~\cite{TRICOCHE2002topologytracking,Post2003SOFlowVis,theisel2003feature,Tino2011StableFeatureFF,Bujack2020TimeFlowSTA}.} 

Our method is inspired by the notion of tracking graphs that represent the correspondences and tracks of individual connected components of spatial features. We formulate the critical point tracking problem as tracking zeros of a time-varying vector field. To achieve this, we use intuition, definitions, and analyses from previous work on tracking critical points in vector fields using robustness~\cite{chazal2011computing,Primoz2015RobustnessSimplification4VF}. 

\begin{figure}
\centering
\includegraphics[width = \linewidth]{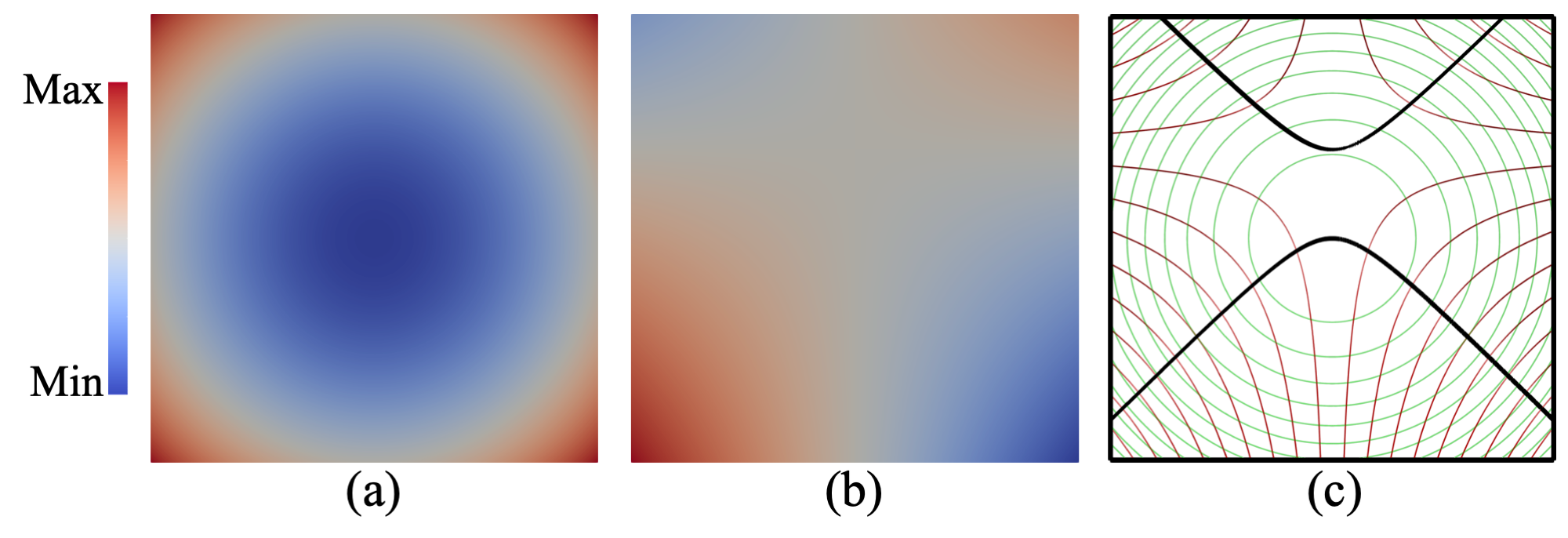}
\caption{Jacobi set of a 2D synthetic bivariate field. (a)~Scalar field $f(x,y) = x^2 + y^2$ (b)~Scalar field $g(x,y) = x(y-8)$ (c)~Level sets of $f$ (green) and $g$ (red). The Jacobi set $\mathbf{J}(f,g)$, shown in black, passes through the points where gradients of $f$ and $g$ align.}
\label{fig:background}
\end{figure}

\vspace{-0.2in}

\section{Background}
\label{sec:background}
This section presents the necessary definitions, notations, and mathematical preliminaries for describing Jacobi set simplification. It introduces critical points of a scalar field and the notion of Jacobi set for a bivariate and a time-varying scalar field, both in the smooth and PL setting~\cite{edelsbrunner2022computational}.

\myparagraph{Scalar field.} 
A \emph{scalar field} $f:\mathbb{M} \rightarrow \mathbb{R}$ is a real-valued function defined on an $n$-dimensional domain $\mathbb{M}$ and maps a point on the domain to a single scalar value. The scalar values usually represent physical quantities such as pressure, temperature, speed, or height, although they may also be synthetically generated.
\autoref{fig:background}(a) shows a 2D scalar field, $f(x,y) = x^2 + y^2$. A point $p \in \mathbb{M}$ is a \emph{critical point} of $f$ if and only if the gradient of $f$ at $p$ vanishes. Else, $p$ is a \emph{regular point}. The scalar value $f(p)$ is called a critical or regular value, respectively. A critical point $p$ is said to be non-degenerate if the Hessian at $p$ is non-singular. If all the critical points of $f$ are non-degenerate and have distinct values, then $f$ is called a \emph{Morse function}. A key result from Morse theory states that all smooth functions are either Morse or can be perturbed using an infinitesimally small quantity into a Morse function. All scalar functions considered in this paper are assumed to be Morse functions.

\myparagraph{Bivariate field.} 
A \emph{multivariate} function or a \emph{multifield} is a collection of scalar functions / scalar fields defined over a common spatial domain. A \emph{bivariate} field $\{f, g\}: \mathbb{M} \rightarrow \mathbb{R}^2$ is a specific case where two scalar fields are defined over a domain $\mathbb{M}$. Bivariate field analysis becomes beneficial when important data features depend on the interaction between two individual fields as opposed to being characterized independently by the two fields. \autoref{fig:background} shows two scalar fields that can either be studied independently or as a bivariate field $\{f,g\}$.

\myparagraph{Jacobi set.} 
Given a bivariate field $\{f, g\} : \mathbb{M} \rightarrow \mathbb{R}^2$ such that the intersection of the sets of the critical points of $f$ and $g$ is a null set, the \emph{Jacobi set} $\mathbb{J} = \mathbb{J}(f, g) = \mathbb{J}(g, f)$ is the set of points where gradients of $f$ and $g$ are linearly dependent.
\begin{equation}
\begin{split}
\mathbb{J}(f,g) = \{x \in \mathbb{M}\ |\ \nabla f + \lambda \nabla g = 0 \text{  or}\\
\lambda \nabla f +  \nabla g = 0 \} 
\end{split}
\end{equation}
for some $\lambda \in \mathbb{R}$. In other words, 

\begin{equation}
\label{eq:def2}
\begin{split}
   \mathbb{J}(f,g) = \{x \in \mathbb{M}\ |\ \text{$x$ is a critical point of } f + \lambda g \text{  or}\\
   \lambda f +  g \} 
\end{split}
\end{equation}

\autoref{eq:def2} essentially characterizes the Jacobi set as the set of critical points of a collection of scalar fields parameterized by $\lambda \in \mathbb{R}$.
Let $k \in \mathbb{R}$ be a regular value of $g$. The scalar field $f_k: g^{-1}(k) \rightarrow \mathbb{R}$ is called the restriction of the scalar field $f$ to the level set $g^{-1}(k)$. The Jacobi set is equivalently defined as
\begin{equation}
\label{eq:def3}
\mathbb{J}(f,g) = cl(\{p \in \mathbb{M}\ |\ p \mbox{ is a critical point of } f_k, k\in\mathbb{R}\})
\end{equation}

The \textit{Smooth Embedding Theorem}~\cite{edelsbrunner2004jacobi} states that, under generic conditions, the Jacobi set of two Morse functions is a smoothly embedded $1$-manifold in $\mathbb{M}$. \autoref{fig:background}(c) shows the level sets of two scalar fields $f$ and $g$ using isocontours. The Jacobi set (black) of the bivariate field $\{f,g\}$ is a curve on the plane. The individual scalar fields are shown in \autoref{fig:background}(a) and \autoref{fig:background}(b).
\begin{figure}
\centering
\includegraphics[width = \linewidth]{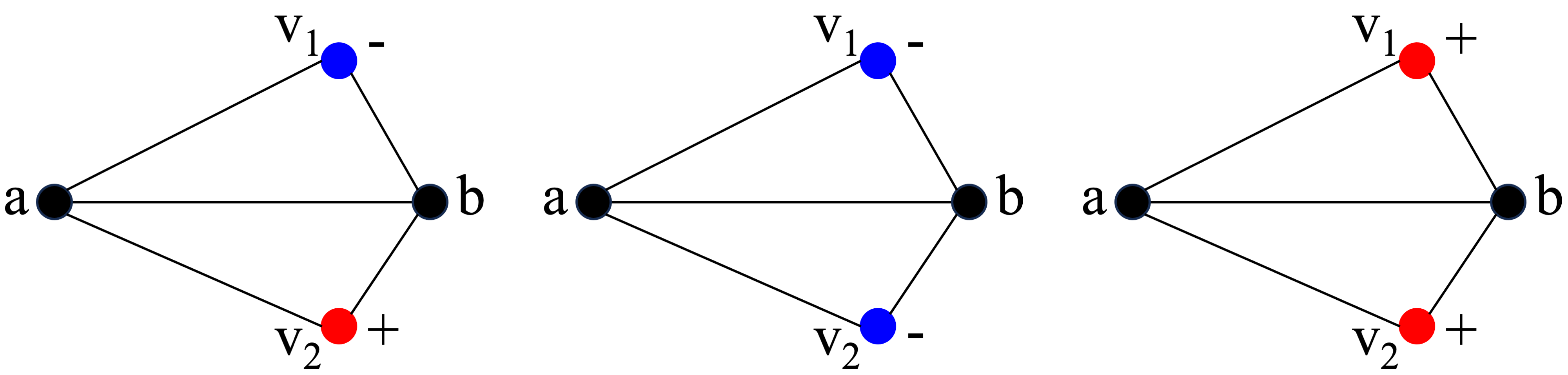}
\caption{Identifying a Jacobi edge. (Left) $h_{\lambda_e}(v_1) < h_{\lambda_e}(a)$ and $h_{\lambda_e}(v_2) > h_{\lambda_e}(a)$, which implies that $ab$ is a regular edge. (Middle)~Both $h_{\lambda_e}(v_1)$ and $h_{\lambda_e}(v_2)$ are smaller than $h_{\lambda_e}(a)$, which implies that $ab$ is a maximum. (Right)~Edge $ab$ is a minimum Jacobi edge.}
\label{fig:background_PL}
\end{figure}

\myparagraph{Piecewise linear function.}
In practice, the scalar fields are often available as a sample over a spatial domain. The spatial domain is represented as a mesh, either a structured grid or a triangulation. The samples of the scalar field are available at the vertices of the mesh. A piecewise linear (PL) approximation of the scalar field is obtained via linear interpolation in the interior of cells of all dimensions. 

Let $f: |K| \rightarrow \mathbb{R}$ be a PL function defined on a triangulation $K$ of a manifold $\mathbb{M}$. Critical vertices of $f$ are characterized by the Betti numbers of its lower link~\cite{delfinado1993incremental,edelsbrunner2022computational}.
Non-Morse functions can be transformed into Morse functions by addressing degeneracies through the simulation of simplicity~\cite{edelsbrunner1990simulation}. In the following, we always assume that the scalar field is a Morse function.

The Jacobi set of two PL functions may be computed as the set of critical points of the family of functions $h_{\lambda} = f + \lambda g$, following \autoref{eq:def2}. According to the \textit{Critical Edge Lemma}~\cite{edelsbrunner2004jacobi}, the Jacobi set lies on the edges of the triangulation $K$. 
If an edge $e=(a,b) \in K$ belongs to the Jacobi set then the function $h_{\lambda}$ is necessarily flat along the edge, \ie $h_{\lambda_e}(a) = h_{\lambda_e}(b)$. This value $\lambda_e$ for which $h_{\lambda}$ is flat can be computed as  
$\lambda_e = (f(b) - f(a)) / (g(a) - g(b))$. Further, the edge $e$ is classified as critical based on a characterization of the lower link of $e$ with respect to $h_{\lambda_e}$. \autoref{fig:background_PL} shows the categorization of an edge $ab$ as regular, maximum, or minimum. If $h_{\lambda_e}(v_1) < h_{\lambda_e}(a)$ and $h_{\lambda_e}(v_2) > h_{\lambda_e}(a)$ then the edge $ab$ is a regular edge. If both $h_{\lambda_e}(v_1)$ and $h_{\lambda_e}(v_2)$ are lower or higher  than $h_{\lambda_e}(a)$, then $ab$ is an extremum edge. The Jacobi set is obtained by stitching the critical edges at their endpoints. The resulting set is a 1-manifold due to the \emph{Even Degree Lemma}, which states that the degree of each vertex in the Jacobi set is even.
\begin{figure}
\centering
\includegraphics[width = \linewidth]{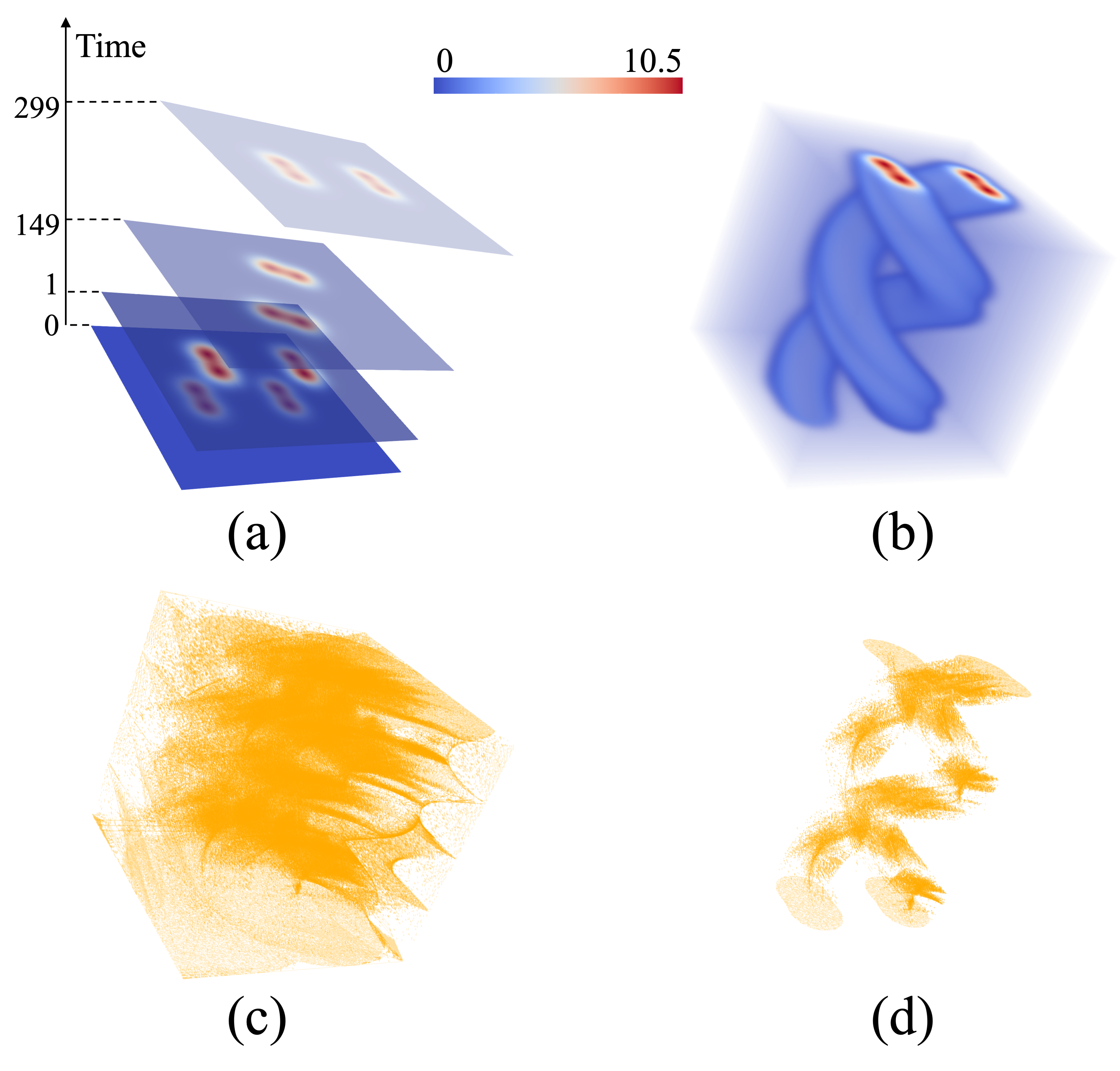}
\caption{Jacobi set of a time-varying scalar field. (a)~A synthetic 2D time-varying scalar field obtained by applying incremental rotations over time. Select time steps within the range [1..299]. (b)~All time steps of the scalar field. stacked together to show the spiral path followed by the two pairs of maxima. (c)~The Jacobi set of the time-varying scalar field is noisy. It is difficult to identify the two primary data features due to clutter. (d)~Edges of the Jacobi set filtered based on the scalar value. Edges with endpoints having scalar value less than $1$ are discarded, leaving only a few Jacobi edges that are in the spatial proximity of important features.}

\label{fig:background_TV}
\end{figure}

\myparagraph{Time-varying scalar field.}
Let $f: \mathbb{M} \times \mathbb{R} \rightarrow \mathbb{R}$ be a time-varying scalar function. The second parameter of $f$ is considered as the time axis. Consider a time function $g: \mathbb{M} \times \mathbb{R} \rightarrow \mathbb{R}$ defined as $g(p,t)  = t$. Following the formulation in \autoref{eq:def3}, the Jacobi set $\mathbb{J}(f,g)$ is the collection of tracks of critical points of $f$ at every time step. This is equal to the set of critical points of $f$ at every time step and the temporal edges that represent the correspondences between critical points lying on two consecutive time steps. 

In this paper, we consider 2D PL time-varying scalar fields $\{f_t\}$, where $f_t$ is a 2D PL scalar field defined in time step $t$. The domain $\mathbb{M} \times \mathbb{R}$ is triangulated by extending a triangulation of the spatial domain. We include temporal edges between corresponding vertices in two consecutive time steps and insert diagonals to decompose the resulting prism into three tetrahedra. The Jacobi set $\mathbb{J} (\{f_t\})$ of the time-varying scalar field is computed over this triangulated domain using the standard PL algorithm~\cite{edelsbrunner2004jacobi}. 

\autoref{fig:background_TV}(a) shows a subset of $300$ time steps of a synthetic time-varying scalar field. The scalar field contains two pairs of Gaussians undergoing rotation over time. In  \autoref{fig:background_TV}(b), all time steps are stacked together, highlighting the rotational pattern formed by these pairs. \autoref{fig:background_TV}(c) shows the Jacobi set, considering time as the second field. However, the Jacobi set is cluttered, and the noisy Jacobi edges occlude the prominent rotational pattern. \autoref{fig:background_TV}(d) shows Jacobi edges close to prominent features by applying a manual threshold. Edges with endpoints having scalar field value less than $1$ are discarded. In this paper, our goal is to simplify the Jacobi set such that the simplified Jacobi set captures the important data patterns.  

\section{Robustness-based simplification of the Jacobi set}
\label{sec:robBasedSimplificaiton}
In this section, we introduce the notion of a simplified Jacobi set for time-varying scalar fields based on the robustness of critical points of a vector field. Critical points of a scalar field occur where the corresponding gradient field vanishes. Extremum critical points of the scalar field correspond to sources and sinks in the gradient field. A saddle corresponds to a saddle in the associated gradient field. Hence, the Jacobi set, consisting of the tracks of critical points of a time-varying scalar field, is the same as the track of zeros of the associated gradient field. We first briefly introduce the notion of robustness of critical points in a vector field and then describe our adaptation to gradient fields followed by a method for Jacobi set simplification.

\begin{figure}
\centering
\includegraphics[width = \linewidth]{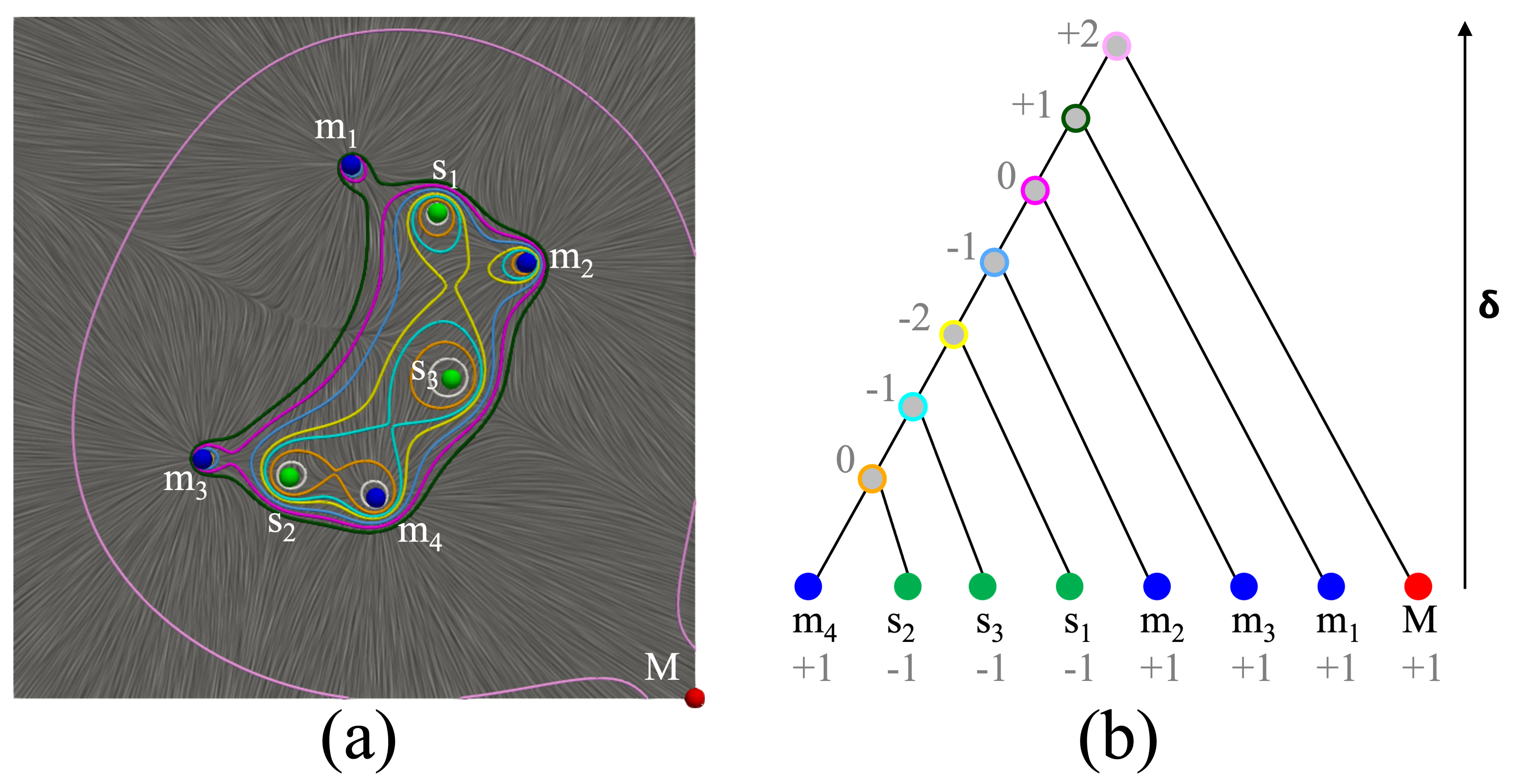}
\caption{(a)~A synthetic 2D vector field visualized using a LIC image and annotated with a subset of its critical points: sources (blue), sink (red), and saddles (green). Isocontours of the vector magnitude field are shown in different colors. The corresponding sublevel sets contain subsets of critical points. All contours corresponding to the sublevel set of a critical point may not be visible. Contours at smaller isovalues may lie close to the critical point and hence are occluded by the glyph used to show the critical point. (b)~The merge tree captures the evolution of the connected components of the sublevel set components. Nodes are annotated with the degree. The degree of the leaf nodes is equal to their Poincar\'{e
} index, +1 for sources/sinks and -1 for saddles. The degree of an internal node is equal to the sum of degrees of all leaf nodes within the subtree rooted at the node.}
\label{fig:robustnessMergeTree}
\end{figure}

\myparagraph{Robustness.} 
In a 2D vector field, the critical point is a point where the vector vanishes and is characterized by its Poincar\'{e} index. Sources and sinks have an index of +1, and a saddle has an index of -1. \autoref{fig:robustnessMergeTree}(a) shows a 2D synthetic vector field. A subset of its critical points, source $M$ (red), sinks $m_i$ (blue), and saddles $s_i$ (green) are highlighted. \emph{Robustness} is a measure that quantifies the stability of a critical point in a vector field. Intuitively, it represents the difficulty or effort required for perturbing away the critical point. 

The formal definition of robustness depends on a topological structure called the merge tree. The \emph{merge tree} of a scalar field $f$ is an abstract representation of the evolution of sublevel set components of $f$. Leaves of the merge tree correspond to the minima of the underlying scalar field and represent the birth of a component. The internal nodes correspond to saddle critical points and represent the merging of sublevel set components.

Consider a 2D vector field $V: \mathbb{R}^2 \rightarrow \mathbb{R}^2$. The robustness of a critical point of $V$ is determined by the critical points that lie within a specific spatial neighborhood of the critical point. This spatial neighborhood is in turn determined by the merge tree of the vector magnitude function $\|V\|$ defined on the plane. A connected component of the sublevel set of $\|V\|$ contains one or more nodes of the merge tree $\mathcal{MT}(\|V\|)$, see \autoref{fig:robustnessMergeTree}. The collection of nodes contained within $C$ together with their incident arcs constitutes a subtree of $\mathcal{MT}(\|V\|)$. Given such a path connected component $C$, the degree $deg(C)$ is defined as the sum of indices of critical points of the vector field $V$ that are enclosed within the component. The degree $deg(C)$ is assigned to the root of this subtree.
Note that critical points of the vector field appear as leaf nodes (minima) of $\mathcal{MT}(\|V\|)$. The \emph{degree} of these leaf nodes is equal to their Poincar\'{e} index. The \emph{degree} of an internal node of $\mathcal{MT}(\|V\|)$ is equal to the sum of degrees of all leaves of the subtree rooted at the node. 

The \emph{static robustness} $\delta$ of a critical point is defined as the vector magnitude of the lowest ancestor node in the merge tree with degree~0~\cite{chazal2011computing}. This paper utilizes only the notion of static robustness and we often refer to it as robustness, without the qualifier. In the vector field shown in \autoref{fig:robustnessMergeTree}, the source $m_4$ has a lower robustness value as compared to the other sources. \newtext{Note that static robustness of zeros in vector fields is theoretically different from the concept of persistence of critical points in scalar fields~\cite{chazal2011computing,edelsbrunner2011quantifying,Primoz2015RobustnessSimplification4VF}. We discuss the difference between the two approaches towards Jacobi set simplification in \autoref{sec:analysis}}.

 %

 \myparagraph{Jacobi set and critical point tracks.} 
 Given a 2D PL time-varying scalar field $\{f_t\}$, where $f_t:\mathbb{R}^2 \rightarrow \mathbb{R}$ is the scalar field in the $t^{th}$ time step, we derive the Jacobi set $\mathbb{J}(\{f_t\})$  by introducing the corresponding time-varying gradient field $\{\nabla f_t\}$, where $\nabla f_t:\mathbb{R}^2 \rightarrow \mathbb{R}^2$. The critical points of $f_t$ are equal to the points where $\nabla f_t$ vanishes. Minima, maxima, and saddles of the scalar field behave as sources, sinks, and saddles of the corresponding gradient field, respectively. These critical points are the 0-value minima of the gradient magnitude field $G_t = \|\nabla f_t\|$. Hence, the set of tracks of critical points of $\{f_t\}$ is the same as the set of tracks of zeros of $\{G_t\}$. We denote the original Jacobi set of the time-varying scalar field either as $\mathbb{J}(\{f_t\})$ or $\mathbb{J}(\{G_t\})$.
 
\myparagraph{Simplifying the Jacobi set.} 
A \emph{$\delta$-sublevel set} of the gradient magnitude field $G_t$ is the preimage $G_t^{-1} (-\infty,\delta]$. 
Let $\mathcal{C}(G_t,\delta)$ denote the set of connected components of the $\delta$-sublevel set. Each connected component in  $\mathcal{C}(G_t,\delta)$ contains a collection of critical points of the corresponding gradient vector field $\nabla f_t$, similar to the illustration in \autoref{fig:robustnessMergeTree} for a generic vector field. These critical points correspond to critical points of the scalar field $f_t$ for all $\delta \geq 0$ and for all time $t$. Consequently, all edges of the Jacobi set $\mathbb{J}(\{f_t\})$ pass through the collection $\mathcal{C}(G_t,\delta)$ for all $t$. The collection of edges of $\mathbb{J}(\{f_t\})$ passing through a pair of spatially overlapping components from adjacent time steps may be represented by a single Jacobi edge, resulting in a reduction in the size of the Jacobi set. \newtext{Figure~1 in the supplementary material illustrates the simplification.} Increasing $\delta$ results in larger components of $\mathcal{C}(G_t,\delta)$, which in turn implies a larger number of critical points clustered together into a single component and represented by a single Jacobi edge. This implies that $\delta$ is a natural parameter for driving Jacobi set simplification.

Given a vector-valued function $V: \mathbb{R}^2 \rightarrow \mathbb{R}^2$, the function $\widetilde{V}: \mathbb{R}^2 \rightarrow \mathbb{R}^2$ is called a \emph{$\delta$-perturbation} of $V$ if $$ \sup_{p \in \mathbb{R}^2} \lVert V(p) - \widetilde{V}(p) \rVert \leq \delta. $$ 
The notion of $\delta$-perturbation aids in the numerical realization of the process that removes critical points with low robustness.  Specifically, it quantifies the effort required to remove the critical points. We transport the notion of robustness to time-varying scalar functions and to their Jacobi set by studying and characterizing $\delta$-robust critical points of the associated gradient vector field. 

A \emph{valid simplification} $\mathbb{J}^*$ of $\mathbb{J}(\{G_t\})$ is a Jacobi set that is derived from $\mathbb{J}(\{G_t\})$ via a $k\delta$-perturbation of $\nabla f_t$ and satisfies the following conditions:
\begin{enumerate}
    \item $\mathbb{J}^* = \mathbb{J}(\{\widetilde{G}_t\})$, where $\widetilde{G}_t$ is the gradient magnitude field of a $k\delta$-perturbation ($k$ is a small positive integer) of $\nabla f_t$ that consists of only $\delta$-robust critical points.
    \item $\mathcal{C}(G_t,\delta) = \mathcal{C}(\widetilde{G}_t,\delta)$  for all $t$.
    \item  $\phi_t(\widetilde{G}_t)(x) \leq \phi_t(G_t)(x)$ for all $t$, where the function $\phi(G_t)$ maps a pair of $\delta$-sublevel set components to the number of Jacobi edges between them, $\phi_t(G_t) : \mathcal{C}(G_t,\delta) \times \mathcal{C}(G_{t+1},\delta) \rightarrow \mathbb{Z}^{+} \cup \{0\}$.
\end{enumerate}
 
The time-varying gradient magnitude field $G_t$ and its perturbation  $\{\widetilde{G}_t\}$ are close to each other, resulting in a Jacobi set $\mathbb{J}^*$ that is close to $\mathbb{J}(\{G_t\})$. The first condition also ensures that noisy critical points, that are not $\delta$-robust, are perturbed away and the stable critical points form the simplified Jacobi set. The second condition ensures that $G_t$ is not altered significantly and the geometry of $\mathbb{J}(\{G_t\})$ is closely approximated. Finally, the third condition ensures a reduction in the number of Jacobi set edges by considering each pair of $\delta$-sublevel set components from adjacent time steps that contributes an edge to the Jacobi set.

\section{Algorithm}
\label{sec:algorithm}
We now present a method to compute a valid simplification of a given Jacobi set $\mathbb{J}(\{G_t\})$. The input scalar field $\{f_t\}$ at every time step is defined on the vertices of a grid and is interpolated within the interior of cells. The algorithm consists of the following three steps:

\myparagraph{Gradient computation.} 
A gradient field is computed for each PL function $f_t$ via central differences. At a grid vertex, the partial derivative along the X-axis is computed as the difference between the scalar values of its left and right neighbors. Similarly, neighbors along the Y-axis are used to compute the partial derivative along the Y-axis. The two partials constitute the gradient vector. The sequence of gradient fields $\{\nabla f_t\}$ and the corresponding gradient magnitude fields $\{G_t\}$ are computed in the first step. For smooth functions, the critical points of a scalar field $f_t$ correspond to a zero of the gradient vector field $\nabla f_t$. However, such a clear correspondence may not exist for a PL scalar field. So, we assign a zero value for $\nabla f_t$ and $G_t$ at vertices that are identified as critical points of $f_t$. $G_t$ is evaluated at all grid vertices and linearly extended in the interior of the cells.

\myparagraph{$\delta$-sublevel set computation.} 
In the second step, we compute the connected components of the $\delta$-sublevel sets of $\{G_t\}$.  \autoref{sec:implemetation} provides the implementation details. 

\myparagraph{Simplified Jacobi set computation.} 
The simplified Jacobi set is computed as the tracking graph of the $\delta$-sublevel set components of $\{G_t\}$. Correspondence between two $\delta$-sublevel set components from adjacent time steps is determined via spatial overlap. 
If $C_1$ and $C_2$ are two components that have a non-empty spatial overlap and further the degrees $deg(C_1) \not = 0$ and $deg(C_2) \not = 0$, then we insert an edge of the Jacobi set between $C_1$ and $C_2$. A geometric embedding of this edge is computed by using the centroid of the component as a representative and inserting the edge connecting the two representatives. If the component $C_1$ from the earlier time step has a non-empty overlap with multiple components from the next time step, then we choose the component with the largest overlap in terms of area.

The Jacobi set is hence computed as a collection of edges between representatives of $\delta$-sublevel set components. The resulting set may contain short or broken tracks, due to short-lived critical points caused by noise. A post-processing step addresses this issue as discussed in \autoref{sec:implemetation}.

\section{Analysis}
\label{sec:analysis}
In this section, we present a theoretical analysis of the algorithm based on a study of the stability of the zeros in the gradient field. We begin by stating two results regarding Riemann maps~\cite{conway2012functions,conway2012functions2}.

\begin{theorem}[Riemann Mapping Theorem~\cite{conway2012functions}]
 Let $A \subset \mathbb{C}$ be a simply connected region which is not the whole plane and let $a \in A$. There exists a unique analytic function $R : A \rightarrow \mathbb{C}$ satisfying the following properties:
 \begin{itemize}
     \item $R(a) = 0$ and the derivative $R'(a)>0$
     \item $R$ is injective
     \item $R(A) = \mathbb{D}$, the unit open disk in the plane
 \end{itemize}
 \end{theorem}

\begin{theorem}[~\cite{conway2012functions2}]
 If $A$ is a bounded simply connected region and $R^{-1}:\mathbb{D}\rightarrow A$ is the Riemann map with $R^{-1}(0) = 0$ and $(R^{-1})'(0)>0$, then $R^{-1}$ extends to a homeomorphism of $\overline{\mathbb{D}}$ onto $\overline{A}$ iff the boundary of $A$, $Bd(A)$, is a Jordan curve. Further, $Bd(A)$ is homeomorphic to $Bd(\mathbb{D})$.
\end{theorem}

Given a connected component of a $\delta$-sublevel set of the gradient magnitude field $G = \|\nabla f\|$ whose degree is non-zero, we will use the Riemann mapping theorem to show the existence of a simplified vector field. This simplified vector field will contain a single zero within each component of the $\delta$-sublevel set $\mathcal{C}(G,\delta)$. We \ prove the following result:


\begin{theorem}
Let $X \in \mathcal{C}(G,\delta)$ be a simply connected component of the $\delta$-sublevel set of the gradient magnitude field G, with non-zero degree and its boundary is a Jordan curve. There exists a continuous vector field $\widetilde{V}:\overline{X} \rightarrow \mathbb{R}^2$ that satisfies the following properties:
\begin{itemize}
    \item $\widetilde{V}$ is a $2\delta$-perturbation of $\nabla f_{|\overline{X}}$,
    \item $\widetilde{V}_{|Bd(X)} = \nabla f_{|Bd(X)}$, and
    \item $\widetilde{V}$ has exactly one singularity.
\end{itemize}
\end{theorem}
\begin{proof}
We first reduce the problem onto a simpler domain, namely the disk, rather than dealing with a generic simply connected Jordan region $X$. The proof consists of two steps.

\myparagraph{Step~1: Reduce the problem to a disk.}
Theorem~1 and~2 guarantee the existence of a homeomorphism $R: \overline{X} \rightarrow \overline{\mathbb{D}}$, the extended Riemann map. Now, define a vector field on the closed unit disk $U: \overline{\mathbb{D}} \rightarrow \mathbb{R}^2$ as a composition $U = \nabla f_{|\overline{X}} \circ R^{-1}$, which allows us to transfer a simplified vector field designed within $\overline{\mathbb{D}}$ onto $\overline{X}$ as shown in \autoref{fig:theoretical-analysis-illustration}.

\myparagraph{Step~2: Construct a simplified vector field on the unit disk.}
We construct a continuous vector field $\widetilde{V}_D: \overline{\mathbb{D}} \rightarrow \mathbb{R}^2$ on the closed disk as $\widetilde{V}_D(r,\theta) = rU(1,\theta)$. This function is composed with $R$ to give $\widetilde{V} = \widetilde{V}_D \circ R : \overline{X} \rightarrow \mathbb{R}^2$, see \autoref{fig:theoretical-analysis-illustration}
\begin{figure}
\centering
\includegraphics[width = \linewidth]{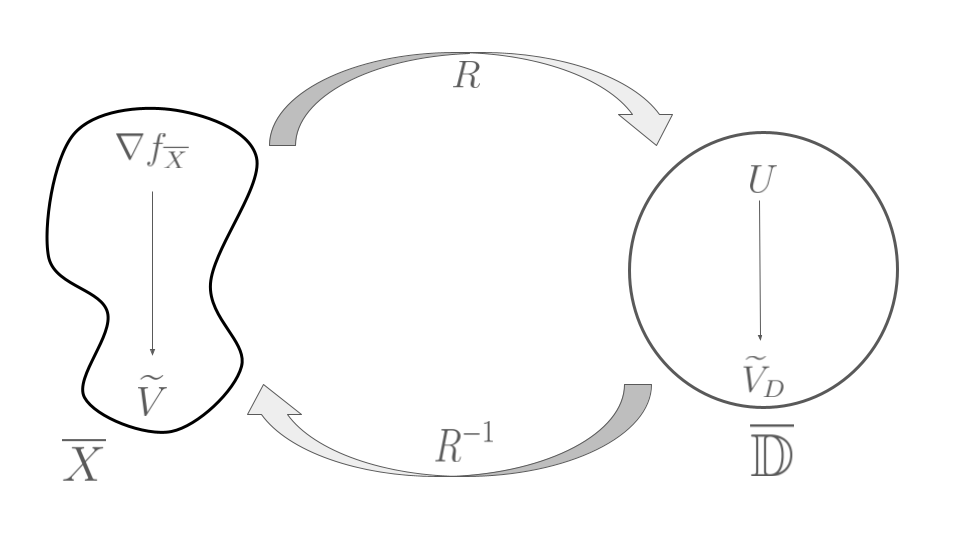}
\caption{A simplified gradient vector field is constructed within the unit disk $\overline{\mathbb{D}}$ and then transferred to the region $\overline{X}$.}
\label{fig:theoretical-analysis-illustration}
\end{figure}
The function $\widetilde{V}$ satisfies the three properties from Theorem~3. Since $U$ has constant magnitude $\delta >0$ on $Bd (X)$ and hence does not have any zeros, it follows that $\widetilde{V}_D$ and $\widetilde{V}$ have exactly one singularity. Further, it is evident from the construction that $\widetilde{V}_{D|Bd(\mathbb{D})} = U_{|Bd(\mathbb{D})}$ and hence $\widetilde{V}$ and $\nabla f_{| \overline{X}}$ are equal on $Bd(X)$. Finally, since $X$ is a $\delta$-sublevel set, the following is true for all values of $0 \leq r \leq 1$ and $0 \leq \theta \leq 2\pi$ 
\begin{equation*}
\|U(1,\theta)\| = \delta \implies \|rU(1,\theta)\| \leq \delta \implies \|\widetilde{V}_D\| \leq \delta 
\end{equation*}
Applying the triangle inequality, we have
\begin{align*}
\|\widetilde{V}_D(r,\theta) - U(r,\theta)\| &\leq 2\delta, \mbox{and} \\
\|\widetilde{V}(x) - \nabla f(x)\| &\leq 2\delta
\end{align*}
Hence, the simplified vector field $\widetilde{V}$ within $X$ is a $2\delta$-perturbation of the gradient field, and its magnitude is smaller than $\delta$.
\end{proof}

Theorem~3 guarantees the existence of a simplified continuous vector field that differs from the gradient of the input scalar field only within the $\delta$-sublevel set components. Further, each zero of the simplified vector field within a $\delta$-sublevel set component is at least $\delta$-robust. 

Next, we consider the case where the connected component of the $\delta$-sublevel set has degree-0. We use an adaptation of a theorem developed in the study of robustness for vector fields~\cite{chazal2011computing} to show that the gradient field can be simplified within such components to remove all critical points. 

\begin{theorem}[~\cite{chazal2011computing}]
Let $X \in \mathcal{C}(G,\delta)$ and let $\nabla f$ have degree-0 in $X$. There exists a continuous vector field $\widetilde{V} : \overline{X} \rightarrow \mathbb{R}^2$ that satisfies the following properties:
\begin{itemize}
    \item $\widetilde{V}$ is a $\delta$-perturbation of $\nabla f_{|\overline{X}}$,
    \item $\widetilde{V}_{|Bd(X)} = \nabla f_{|Bd(X)}$, and
    \item $\widetilde{V}$ has no singularity.
\end{itemize}
\end{theorem}

Theorem~4 implies that perturbing $\nabla f$ within all such regions $X$ results in a simplified continuous vector field that preserves the $\delta$-sublevel sets. The interpretation of static robustness of a critical point as its stability or the perturbation required to remove it also follows from Theorem~4.

The above results guarantee that the zeros of $\nabla f$ can be perturbed away to retain either a single stable zero (Theorem~3) or remove all zeros (Theorem~4). The perturbation results imply that the zeros (critical points of $f$) within the $\delta$-sublevel set are less significant, need not be tracked individually, and can be replaced by a single stable zero. This corollary drives our simplification algorithm. 
The algorithm tracks connected components of $\mathcal{C}(G_t,\delta)$ across time $t$ if their degree is non-zero. A single edge replaces the collection of Jacobi edges that are incident on the critical points lying within these $\delta$-sublevel set components. \newtext{The simplified Jacobi set computed by our algorithm adheres to the valid simplification constraints outlined in \autoref{sec:robBasedSimplificaiton} and represents the tracks of zeros of the perturbed vector fields described in Theorems~3 and~4. However, the perturbed fields are not computed explicitly.}

The algorithm reduces the number of edges and hence computes a simplified Jacobi set. However, it captures the macroscopic trends of movement of the critical points. Birth and death of less significant critical points do not affect the result and do not appear in the simplified Jacobi set because they occur in the interior of $\delta$-sublevel sets.
\newtext{The $\delta$-sublevel sets provide a clustering of critical points that are in geometric proximity to each other and evolve in a coherent manner over time. Thus, the algorithm ensures temporal and spatial consistency of the simplified Jacobi set while reducing the number of tracks. In contrast, tracking methods that rely on standard persistence-based filtering do not require persistence pairs to be consistently paired over time, nor do the pairs have to be in geometric proximity to each other. Also, the persistence values of critical points may vary significantly over time.}

\section{Implementation}
\label{sec:implemetation}
This section describes details of an implementation of the Jacobi set simplification algorithm,  optimizations that improve the runtime performance and the visual clarity of the simplified Jacobi set. All visualizations in this paper are generated using Paraview\cite{Ayachit2015paraview} and \textsc{TTK}~\cite{Tierny2018ttk}. In the following discussion, $T$ denotes the total number of time steps and $n$ is the total number of vertices in the input 2D grid. 

\myparagraph{Persistence simplification.}
Simulation of simplicity, a symbolic perturbation of the scalar field, is used to handle degeneracies when classifying critical points. However, this results in the creation of spurious critical points within regions where the scalar function is flat. These critical points of $f_t$ are removed within each time step using a persistence-driven topological simplification~\cite{ELZ02}. A small persistence threshold is sufficient to remove these critical points. This helps reduce unnecessary processing of these spurious critical points in subsequent steps. 

\myparagraph{Critical point computation.}
The critical points and their types (maxima, minima, or saddle) are computed and stored. This requires $O(nT)$ time. Storing the type of critical point helps optionally filter the Jacobi set to show edges corresponding to a specific type.

\myparagraph{Gradient field computation.}
Computing the gradient magnitude field $G_t$ at each time step $t$ and explicitly setting the gradient magnitude of critical points to zero requires a runtime of $O(nT)$.
\begin{figure}
\centering
\includegraphics[width = 0.8\linewidth]{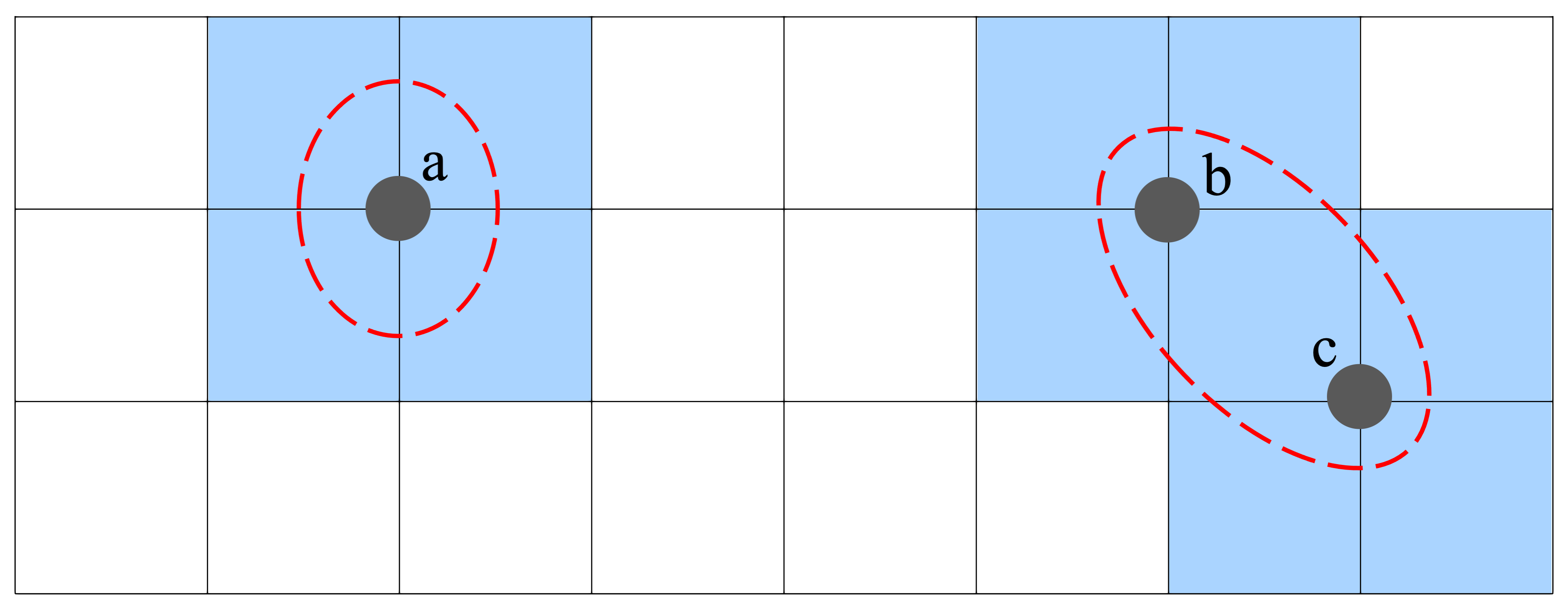}
\caption{$\delta$-sublevel set computation. The $\delta$-isocontour (red) of the gradient magnitude field has two components and the $\delta$-sublevel set contains three critical points $a, b$, and $c$. The algorithm returns the collection of blue cells, a conservative approximation of the $\delta$-sublevel set.}
\label{fig:deltaSublevelset}
\end{figure}

\myparagraph{$\delta$-sublevel set computation.} 
We compute the $\delta$-sublevel sets at each time step $t$ in a single pass over the list of all vertices. If $G_t(p) \leq \delta$ for a vertex $p$ then all cells containing $p$ are inserted into the reported $\delta$-sublevel set. This reported $\delta$-sublevel set is a superset of the true $\delta$-sublevel set. The additional cells ensure that a connected component consisting of a single point is also represented by a spatial region and may be tracked via spatial overlap. \autoref{fig:deltaSublevelset} illustrates the computation with a small example consisting of two components and three critical points. The running time for this step is $O(nT)$. 

\myparagraph{Connected component computation.}
Connected components of the $\delta$-sublevel sets are computed using the connectivity filter of \textsc{TTK}. Further, all degree-0 components are removed. The filter uses a flood fill algorithm and has a runtime of $O(nT)$. Removing all components with degree 0 also takes $O(nT)$ time. 

\myparagraph{Simplified Jacobi set.}
The spatial overlap based tracking filter in \textsc{TTK}, called ``Tracking from Overlap''~\cite{lukasczyk2017nested}, is used to compute a tracking graph. Edges of the simplified Jacobi set are essentially the edges of this tracking graph. This computation has a worst case runtime of $O(nT)$. \newtext{Further details are available in Section~1 of the supplementary material.}

\begin{figure}
\centering
\includegraphics[width = \linewidth]{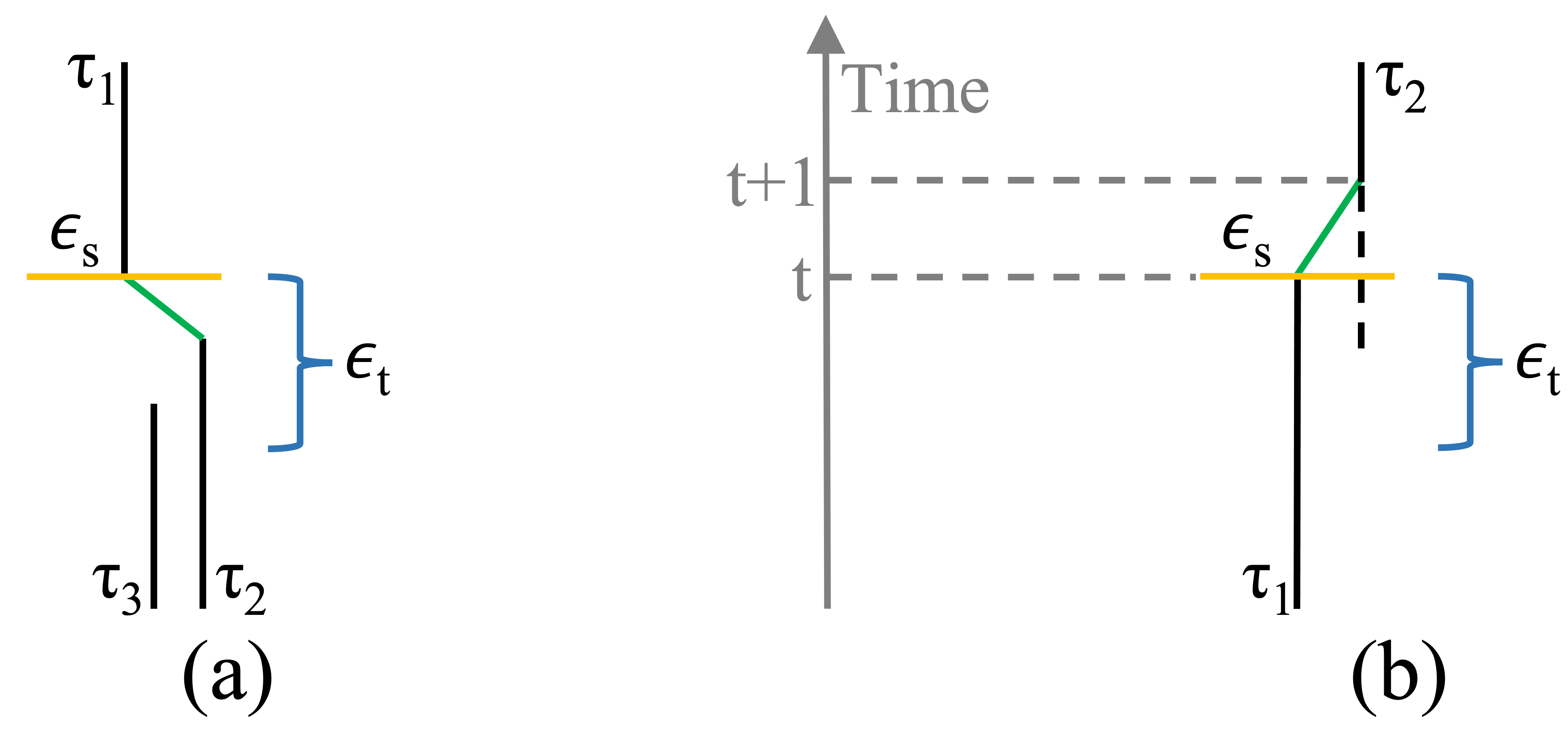}
\caption{Connecting broken tracks based on spatial and temporal proximity. (a)~Tracks $\tau_2$ and $\tau_3$ end within a time interval $\epsilon_t$ of the start of track $\tau_1$ and their end points lie within an $\epsilon_s$-ball (orange) centered at the start point of $\tau_1$. The temporally closer track $\tau_2$ is merged with $\tau_1$ (green edge). (b)~End time of track $\tau_1$ start time of $\tau_2$ are closer than the spatial ($\epsilon_s$) and temporal ($\epsilon_t$) thresholds. They are merged (green edge) and the initial segment of $\tau_2$ (dashed) is discarded.}
\label{fig:postprocessing}
\end{figure}

\myparagraph{Post-processing.}
Noise in the data could result in short-lived or broken tracks. We identify and resolve such issues in a post-processing step to produce the simplified Jacobi set. This step analyzes the proximity of the start and end points of different tracks and determines if they may be merged. Specifically, it checks two criteria as illustrated in \autoref{fig:postprocessing}.
First, for a track $\tau_1$ that begins at time $t$ and at spatial location $(x,y)$, a search within the time interval $[t-\epsilon_t,t]$ determines if there exists another track whose end point lies within an $\epsilon_s$-neighborhood of $(x,y)$ \ie within a ball of radius $\epsilon_s$ centered at $(x,y)$. If such a track $\tau_2$ is found, its endpoint is connected to the start of $\tau_1$, as shown in \autoref{fig:postprocessing}(a). If multiple such tracks are found, the one with the closest end time is declared as $\tau_2$. 
Second, for a track $\tau_1$ that ends at time $t$ and at spatial location $(x,y)$, a search is performed to determine if there exists another track that ends later than time $t$ and whose start time is within the interval $[t-\epsilon_t,t]$ and start point lies within an $\epsilon_s$-neighborhood of $(x,y)$. If such a track $\tau_2$ is found, the end point at time step $t$ in $\tau_1$ is connected to the start point at time step $t+1$ in $\tau_2$ and the segment of $\tau_2$ up to time $t+1$ is discarded. \autoref{fig:postprocessing} (b) illustrates this case.
Significant critical points are expected to correspond to long tracks. Based on this intuition, a length threshold $l$ is applied to discard all short tracks.

In the above implementation, a total of 5 parameters may be tuned to achieve an appropriately simplified Jacobi set. A persistence threshold $\epsilon_p$, the parameter $\delta$, the maximum length $\epsilon_t$ of recoverable temporal breaks in tracks, the size $\epsilon_s$ of the spatial neighborhood used to merge tracks, and a track length threshold $\epsilon_l$. Table 1 in the supplementary material lists the parameter values used in our experiments. Combining the running time of the individual steps, the total runtime of the algorithm is $O(nT)$.

\section{Experimental results}
In this section, we present results on four datasets -- a synthetic rotating Gaussians dataset, 2D von K\'{a}rm\'{a}n vortex street, Boussinesq flow, and sea surface height. We provide qualitative insights into the reduction of the number of Jacobi edges and showcase the benefit of identifying prominent patterns and tracking features corresponding to critical regions. Table 1 in the supplementary material lists the statistics on the size of the simplified Jacobi set, demonstrating the effectiveness and utility of the method in practice. It also lists the parameter values used in the experiments for completeness. 

The simplified Jacobi sets reported here and in subsequent results include tracks of all critical point types -- maxima, minima, and saddles. The implementation may be extended to filter and selectively display tracks that correspond to a specific critical type. All tracks are rendered using cubic spline interpolation. The focus of the prototype implementation in Python is functionality and not efficiency. 
Code optimization using efficient libraries is likely to result in significant runtime improvement. 

\begin{figure}[h!tbp]
\centering
  \includegraphics[width=.9\linewidth]{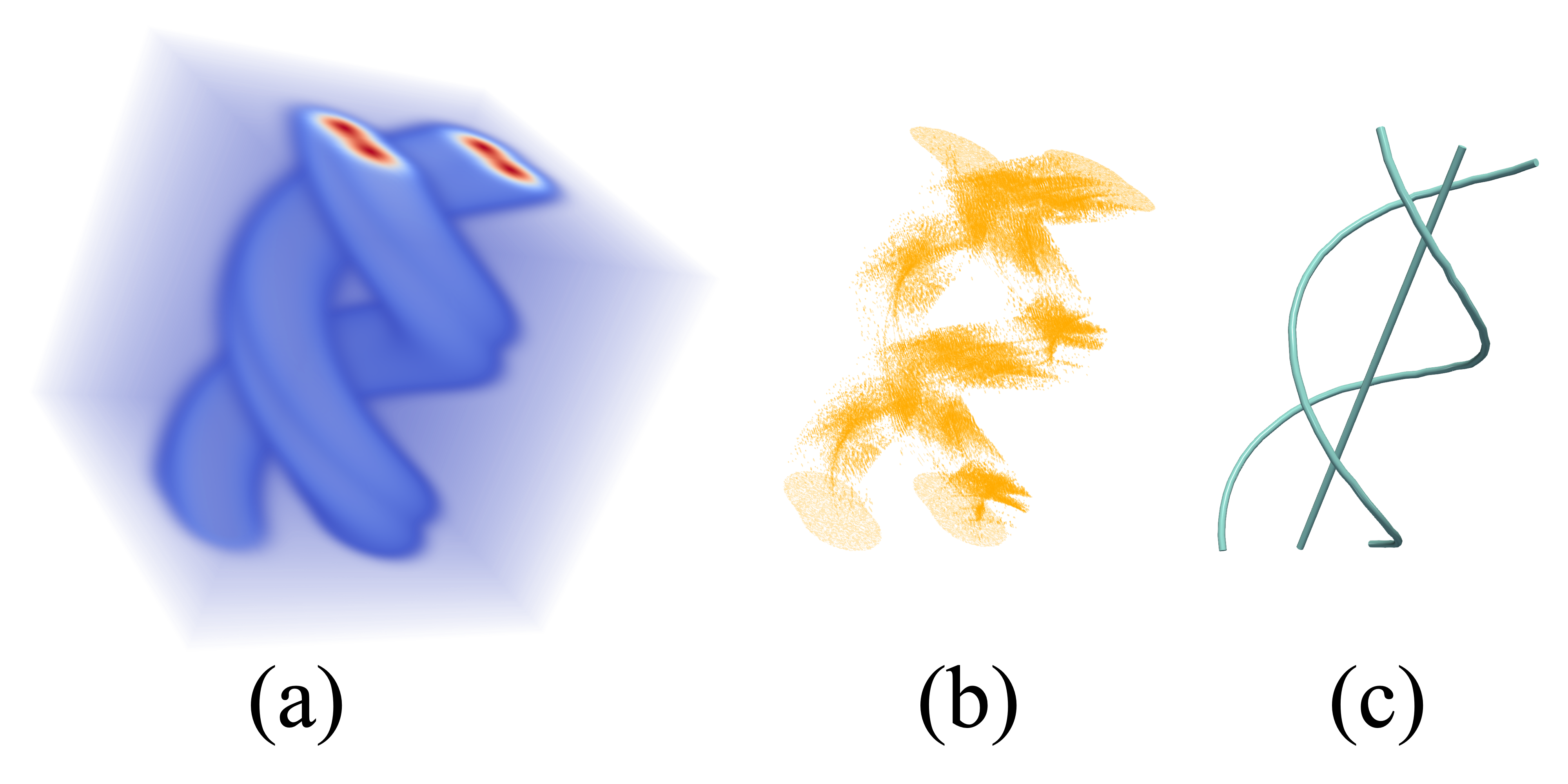}
\caption{(a)~A time-varying scalar field represented as a stack of color-mapped images. The vertical direction corresponds to time. (b)~Edges of the original Jacobi set (orange) after filtering out edges where the scalar value is smaller than 1. (c) The simplified Jacobi set consists of three components. Two tracks capture the spiral pattern of rotating Gaussians. The straight track corresponds to the global minimum.}
\label{fig:Synthetic}
\end{figure}

\myparagraph{Rotating Gaussians.}
A synthetic dataset generated as a time-evolving sum of Gaussians helps analyze the correctness of the algorithm and its implementation. Two pairs of maxima, each consisting of two Gaussian centers (maxima) rotate about the origin over time as shown in \autoref{fig:Synthetic}. The Jacobi set is noisy and cluttered as discussed earlier. We observe that upon filtering the Jacobi set (orange) to retain edges with scalar value greater than $1$, we obtain the spiral pattern followed by the rotating Gaussians. Tracks in the simplified Jacobi set shown in \autoref{fig:Synthetic}(c) effectively capture this pattern, indicating that the simplified Jacobi set is a geometrically accurate representation of the feature (Gaussian center) evolution. Furthermore, a single track for each critical region shows the successful clustering of critical points. The straight track in the center corresponds to the global minimum. No post-processing was required to produce the results.

\begin{figure}[h!tbp]
\centering
  \includegraphics[width=\linewidth]{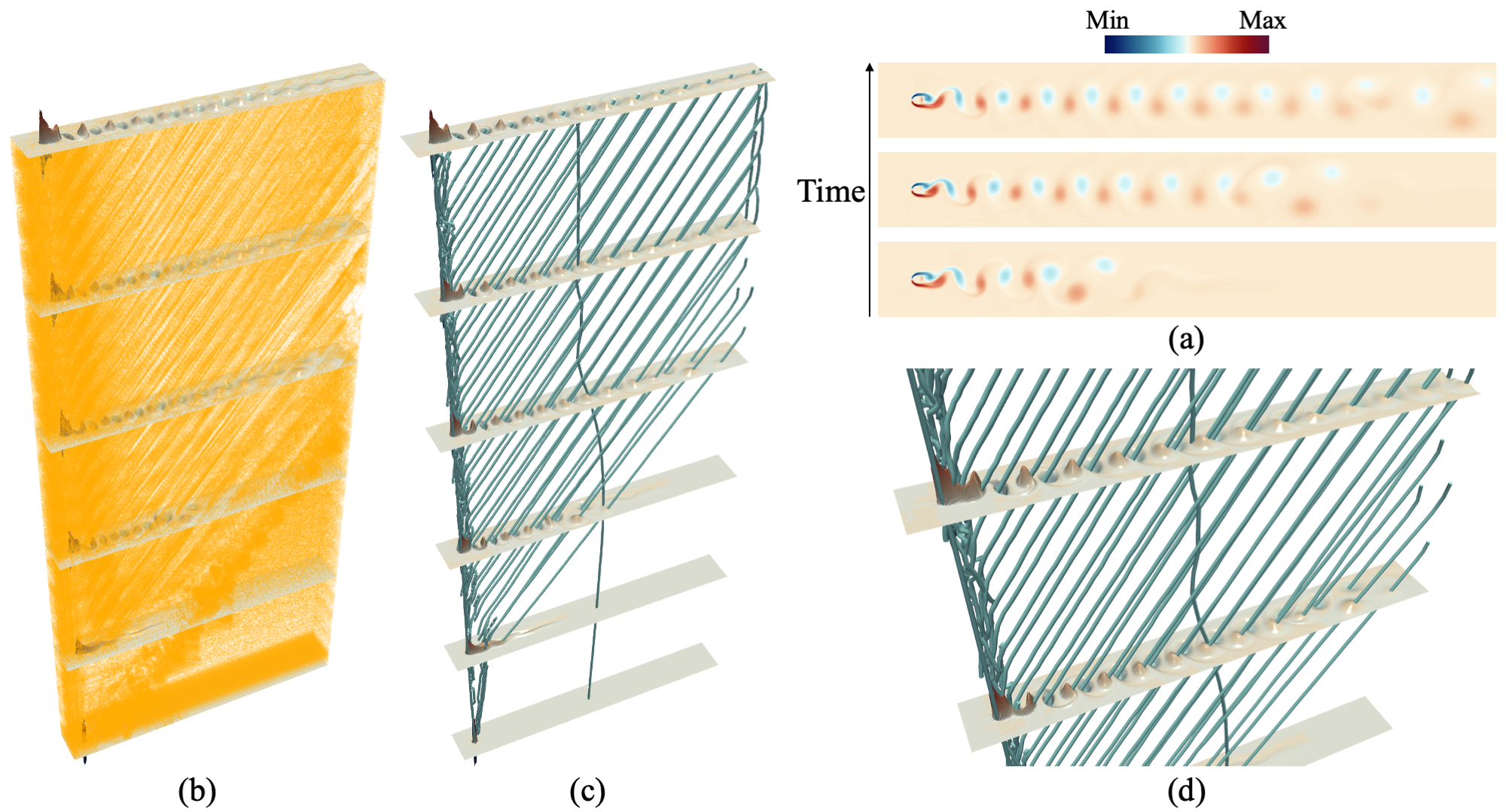}
\caption{(a)~Three time steps from the 2D von K\'{a}rm\'{a}n vortex street dataset color-mapped with $z$-component of vorticity. (b)~Jacobi set is large and noisy. (c)~Simplified Jacobi set is clutter-free. The tracks capture the flow pattern of vortices from left to right.. (d)~A zoomed-in view shows that each track represents an individual critical region around a maximum or minimum, namely a vortex. The vertical track represents a critical region near the domain boundary and corresponds to a saddle.}
\label{fig:Vortex}
\end{figure}
\myparagraph{2D von K\'{a}rm\'{a}n vortex street.}
This dataset of the flow of a fluid around a cylindrical obstacle is generated using Gerris flow solver~\cite{gerrisflowsolver,Guenther17}. The vorticity field is computed using the 2D velocity field available at each grid vertex over 1501 time steps. Since the velocity is a 2D field, only the $z$-component of vorticity is non-zero and is used for the analysis. The dataset is available as a single 3D grid where the $z$-direction denotes time. 

\autoref{fig:Vortex}(a) shows three selected time steps of the vorticity field. A persistence-driven simplification with a small threshold of 0.12\% successfully removes spurious critical points from each time step, but the Jacobi set remains cluttered as shown in \autoref{fig:Vortex}(b). \newtext{A persistence-directed simplification of individual scalar fields may not necessarily produce good results.} The simplified Jacobi set \autoref{fig:Vortex}(c) clearly captures the movement of vortices from left to right. Similar patterns are reported in earlier work, in particular using the lifted Wasserstein matcher~\cite[Figure~14]{Soler2018LiftedWasserstein}. \autoref{fig:Vortex}(d) provides a zoomed-in view showing tracks passing through critical regions of maxima and minima, implying that the geometry of the key tracks is well preserved. The simplified Jacobi set is a clean visual representation with reduced zigzag and minimal clutter, thanks to a substantial reduction (400$\times$) in size. The vertical track passing through the center of the domain represents a cluster of saddles that constitute a large critical region near the domain boundary.

\begin{figure*}[h!tbp]
\centering
  \includegraphics[width=\linewidth]{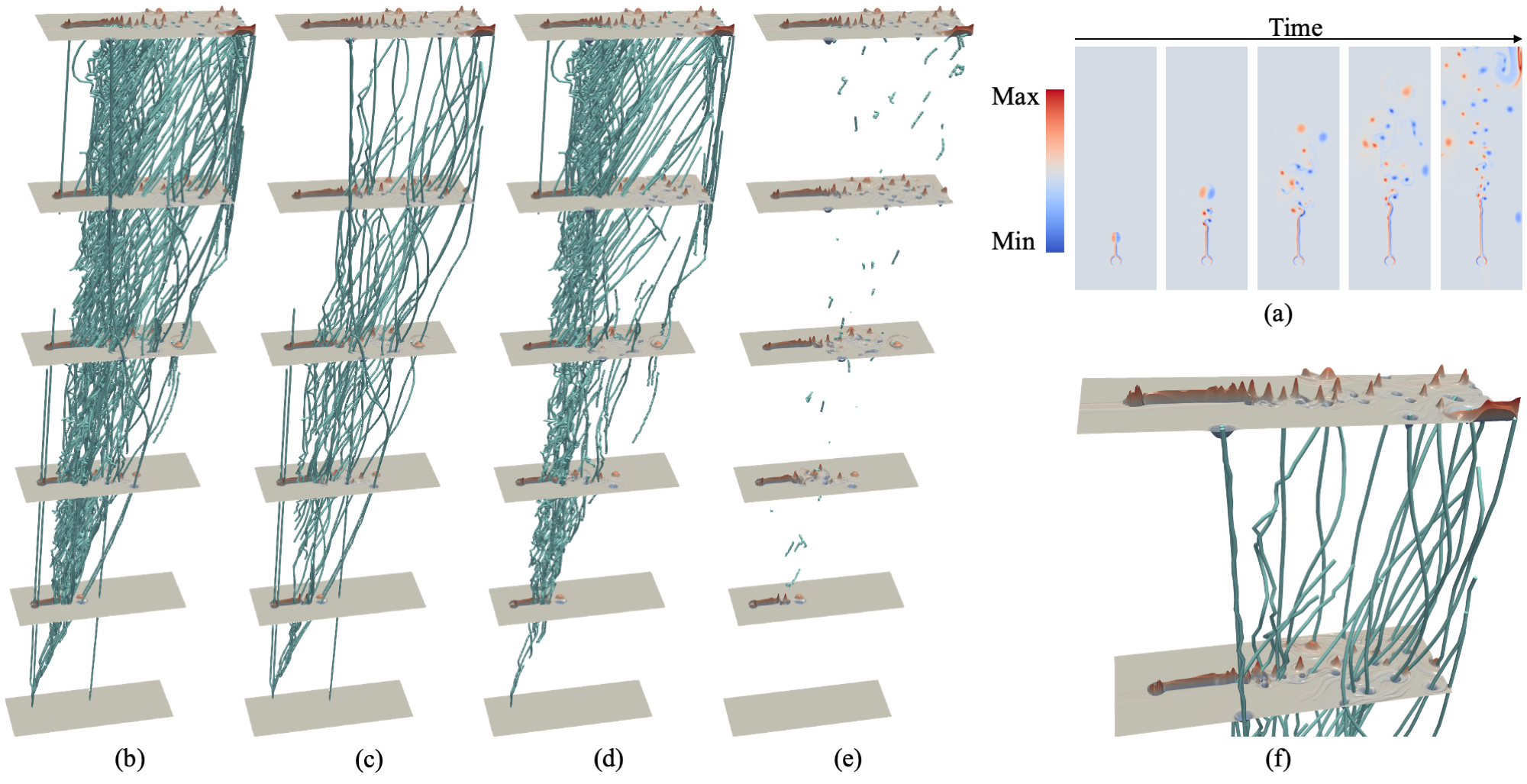}
\caption{(a)~Select time steps of the Boussinesq flow dataset showing the $z$-component of vorticity. (b)~Tracks of all critical regions. (c-e)~Tracks filtered based on their length ($>750, 50-750, <50$). (f)~A close-up view of the long tracks shows that several critical regions are tracked despite the turbulent behavior.}
\label{fig:boussinesq}
\end{figure*}
\myparagraph{Boussinesq flow.}
This dataset representing 2D flow due to a heated cylinder is generated by a simulation using the Gerris flow solver~\cite{gerrisflowsolver,Guenther17}. The vorticity field is computed from the velocity field, as described above, over 2001 time steps and is represented as a scalar field on a 3D grid. \autoref{fig:boussinesq}(a) shows the field in select time steps. 

The simplified Jacobi set represents tracks of critical points of all types. The tracks are classified based on their length and shown in \autoref{fig:boussinesq}(b-e). A visual inspection indicates a significantly larger number of long tracks identified by our algorithm compared to the ones reported by the lifted Wasserstein matcher~\cite[Figure~12]{Soler2018LiftedWasserstein}\newtext{, a method that relies on persistence pair matching for tracking.} Despite the innately unsteady nature of this dataset, our algorithm successfully identifies 38 long tracks that last longer than 750 time steps. \autoref{fig:boussinesq}(f) presents a zoomed-in view of the long tracks passing through important extrema regions, emphasizing the validity of the computed tracks. The number of medium length tracks shown in \autoref{fig:boussinesq}(d), indicates that many critical regions have a life span between 50 and 750 time steps. Experimental results on this turbulent dataset demonstrate the utility of the post-processing heuristic. Tracks in proximity are effectively merged, resulting in a 15$\times$ reduction in the total number of tracks, from 4597 to 312. The number of short tracks is small.

\begin{figure}[h!tbp]
\centering
  \includegraphics[width=.9\linewidth]{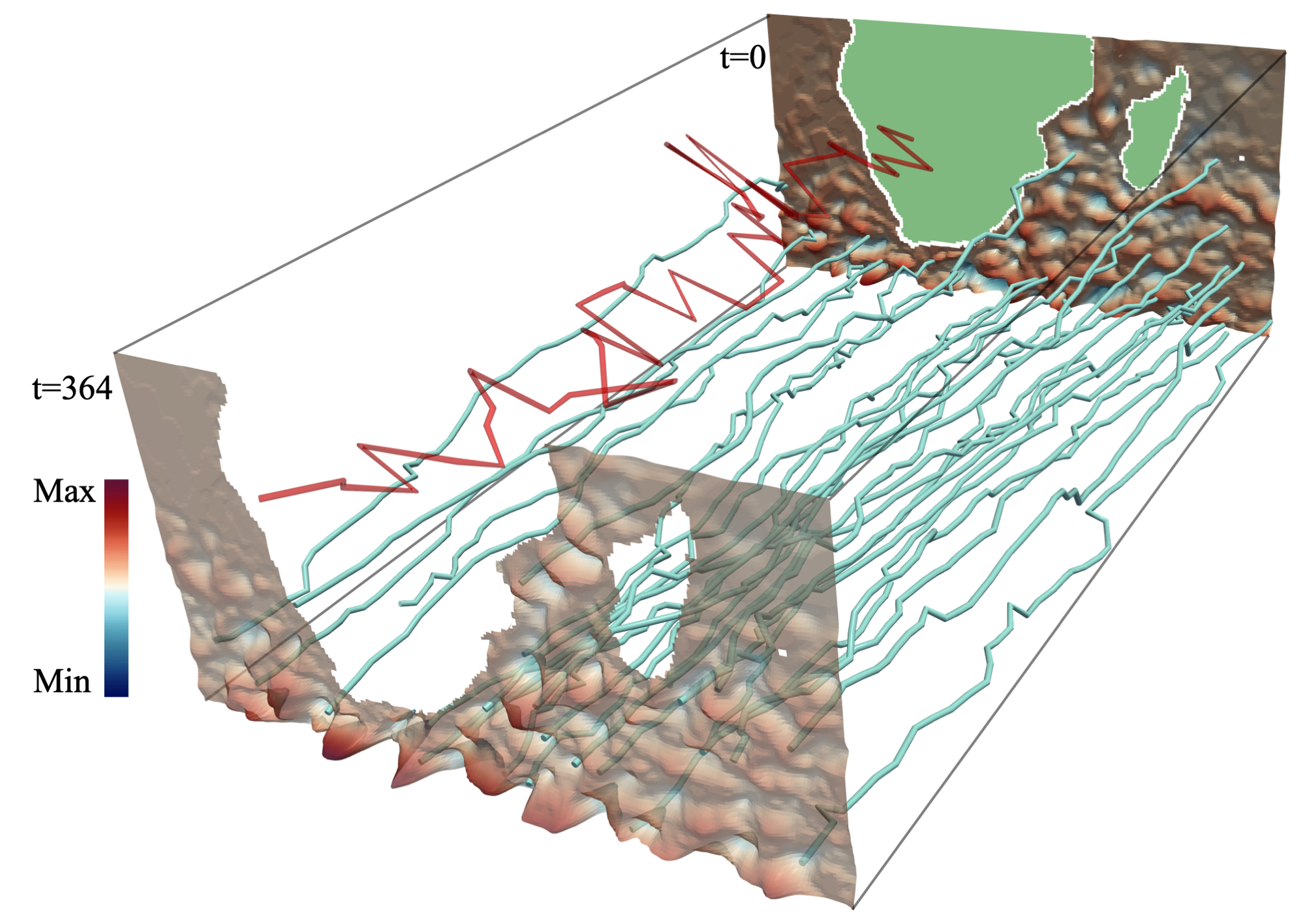}
  \caption{Long tracks ($>$300 days) in the sea surface height dataset. Extrema tracks exhibit a westward drift due to wind. The red track depicts unstable saddles near the coastline.}
\label{fig:sea-surface-height}
\end{figure}
\myparagraph{Sea surface height.}
This dataset from Marine Copernicus~\cite{copernicus} measures the sea surface height from Aug~2021 to Aug~2022. The analysis aims to track significant extrema of the sea surface height over one year to identify features that persist for a long period of time. A length threshold $\epsilon_l = 300$ days is chosen to identify such long tracks. The tracks of extrema shown in \autoref{fig:sea-surface-height} capture the westward drift. This drift is possibly due to the direction of the wind. The simplified Jacobi set is able to capture additional tracks when compared to those extracted using LWM, a direct tracking of critical points~\cite[Figure~11]{Soler2018LiftedWasserstein}. \newtext{We use data from a different year due to unavailability of data for the year presented in LWM. Nevertheless, it is evident via a visual comparison that LWM misses multiple prominent features.} The noisy red track corresponds to the inconsistent and unstable saddles near the boundary or the coastline. Our algorithm successfully clusters and represents all such saddles with a single track.

\vspace{-0.225in}

\section{Conclusions}
We described a novel Jacobi set simplification method in the context of time-varying scalar fields, where the Jacobi set represents tracks of critical points. The simplification is based on an analysis of the associated gradient field and an adaptation of the notion of robustness for vector fields. Experiments show that the method produces clutter-free tracks and can potentially identify a larger number of feature tracks when compared with a method for tracking critical points. We also presented a theoretical guarantee of the existence of a simplified vector field to support the algorithm. A generalization of the theoretical analysis to 3D and to regions that are not Jordan is non-trivial, because there is no direct extension of the Riemann mapping theorem. 

\begin{acknowledgements}
This work is partially supported by a grant from SERB, Govt. of India (CRG/2021/005278). VN acknowledges support from the Alexander von Humboldt Foundation, and Berlin MATH+ under the Visiting Scholar program. Part of this work was completed when VN was a guest Professor at the Zuse Institute Berlin.
\end{acknowledgements}
\vspace{-0.225in}

\bibliographystyle{spmpsci}
\bibliography{references.bib} 

\appendix



\section*{Supplementary material}

# Supplementary Material for "Jacobi Set Simplification for Tracking Topological Features in Time-Varying Scalar Fields"

Dhruv Meduri, Mohit Sharma and Vijay Natarajan

**Abstract**

This document provides additional implementation details, a table listing parameter values used in the implementation, and statistics on the Jacobi sets pre- and post-simplification for different datasets in the paper "Jacobi Set Simplification for Tracking Topological Features in Time-Varying Scalar Fields".

✦

## 1 SIMPLIFIED JACOBI SET

Figure 1(a) shows a few critical points within the gradient magnitude fields $G_t$ and $G_{t+1}$ in two consecutive time steps. The critical points are clustered into $\delta$-sublevel set components $C_1$ and $C_2$. As mentioned in Section 7 of the paper, the spatial overlap based tracking filter in TTK, called "Tracking from Overlap" [1] is utilized to identify overlapping components in the adjacent time steps. The critical points enclosed in the individual components are replaced by their centroid, serving as a single node in the simplified Jacobi set. The Jacobi edges passing through overlapping $\delta$-sublevel set components in adjacent timesteps are replaced by a single arc passing through the centroids of the sublevel set components. Figure 1(b) shows arcs of the simplified Jacobi set resulting from the above-mentioned procedure.

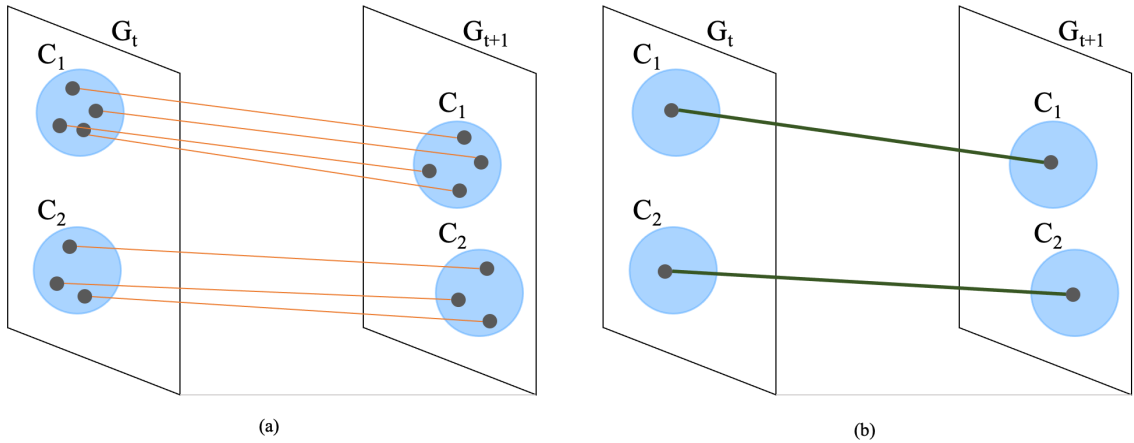

(a)       (b)

Fig. 1: Jacobi set simplification. (a) Two components of $\delta$-sublevel set, $C_1$ and $C_2$, for gradient magnitude field $G_t$ and the $\delta$-sublevel set components in $G_{(t+1)}$ that have a spatial overlap with the former. Jacobi edges (orange) connect critical points lying within $C_1$ and $C_2$. (b) Simplified Jacobi set computed by replacing all orange edges with a single edge (green).

## 2 TRACK STATISTICS

TABLE 1: Track statistics and parameter thresholds. The number of edges of the Jacobi set reduces significantly post simplification, resulting in a clean and clutter-free visualization. The number of tracks #JT-SP in Boussinesq is 38 (large) + 214 (medium) + 60 (short). Parameter thresholds are chosen conservatively, higher values would also work. (#JE: Jacobi edges, #JE-S: Jacobi edges after simplification, #JT-S: Jacobi tracks after simplification, #JT-SP: Jacobi tracks after post-processing)

| Dataset | #JE | #JES | #JT-S | #JT-SP | $\varepsilon_p(\%)$ | $\delta(\%)$ | $\varepsilon_t$ | $\varepsilon_s$ | $\varepsilon_l$ |
|---|---:|---:|---:|---:|---:|---:|---:|---:|---:|
| Rotating Gaussians | 2283031 | 897 | 3 | - | 0.95 | 23 | - | - | - |
| von Kármán vortex street | 13242297 | 33206 | 610 | 57 | 0.12 | 0.13 | 50 | 11 | 150 |
| Boussinesq | 41912203 | 73101 | 4597 | 312 | 7.00 | 0.75 | 200 | 11 | {50,750} |
| Sea surface height | 1051267 | 8894 | 4624 | 27 | 0.20 | 0.02 | 25 | 7 | 300 |



\end{document}